\documentclass[12pt]{article}

\usepackage[T1]{fontenc}
\usepackage[utf8]{inputenc}
\usepackage{authblk}

\usepackage[numbers]{natbib}
\usepackage[OT4]{fontenc}
\usepackage[utf8]{inputenc}

\newcommand{\Exp}{{\rm I\hspace{-0.8mm}E}}
\newcommand{\Prob}{{\rm I\hspace{-0.8mm}P}}

\newcommand{\iz}{{\rm \rlap Z\kern 2.2pt Z}}

\newcommand{\RL}{{\rm I\hspace{-0.8mm}R}}

\newcommand{\proof}{\noindent {\bf Proof:} \ }

\newcommand{\ind}{{1\hspace{-1mm}{\rm I}}}

\newtheorem{theorem}{Theorem}

\newtheorem{proposition}{Proposition}
\newtheorem{definition}{Definition}
\newtheorem{corollary}{Corollary}[section]

\newtheorem{remark}{Remark}
\newtheorem{example}{Example}

\author[1]{Zbigniew Michna}
\author[1]{Wojciech Bombała}
\author[2]{Peter Nielsen}
\affil[1]{Department of Mathematics and Cybernetics

Wrocław University of Economics}
\affil[2]{Department of Mechanical and Manufacturing Engineering Aalborg University}

\title{\bf\LARGE {L\'evy processes in storage\\
and inventory problems}}
\date{}

\begin{document}

\maketitle

\bibliographystyle{abbrv}

\begin{abstract}
In this paper we consider storage and inventory systems. Our aim is to apply and review main
results of the fluctuation theory of stochastic processes in the context of storage and inventory modeling. We describe systems where the inflow is due to a L\'evy process and the outflow is linear and conversely systems where 
the inflow is linear and the outflow is due to a L\'evy process. For such systems we investigate the process of a storage (inventory) level. We give formulas for the probability that the storage level exceeds a certain value. This probability is crucial in the management of storage and inventory systems.

\vspace{5mm}
{\it Keywords: L\'evy process, storage system, inventory system, storage level, probability of an overflow, scale function, gamma process, stable L\'evy process, inverse Gaussian L\'evy process, tempered stable L\'evy process}
\newline
\vspace{2cm}
MSC(2010): Primary 60K25; Secondary 60G51, 60G70.
\end{abstract}

\section{Introduction}
L\'evy processes appear in many theoretical  and practical fields where they serve as a basic skeleton for a description of certain phenomena. They are applied in physics, economics, finance, insurance, queueing systems and other branches. Their features like independence and stationarity of increments or self-similarity in certain cases for instance permit to apply them to model returns of stock prices, claims to insurance companies or an inflow (outflow) to the buffer in queueing (telecommunications) systems. Moreover L\'evy processes serve as a starting point for more complicated models e.g. based on stochastic differential equations. The topic of this research paper is L\'evy processes in the context of inventory management.  

Let us recall the notion of a L\'evy process.
\begin{definition}
A stochastic process $\{X(t), 0\leq t<\infty\}$ with values in $\RL^d$ is a L\'evy process 
if
\begin{itemize}
	\item $X(0)=0$ a.s., 
	\item $X$ has independent increments,
	\item $X(t)-X(s)$ has the same distribution as $X(t-s)$, that is $X$ has stationary increments, 
	\item $X$ is stochastically continuous.
\end{itemize}
\end{definition}

The structure of L\'evy processes is not much complicated.  
L\'evy-It\^o representation shows their stochastic construction which is the following (see e.g. Sato \cite{sa:99})
$$
X(t)=B(t)+\int_{|x|<1}x\,(N_t(dx)-tQ(dx))+\int_{|x|\geq 1}x\,N_t(dx)+at\,,
$$
where $B(t)$ is a Wiener process, $N$ is a point process generated by the jumps of $X$ that is $N=\sum_{\{t:\Delta X(t)\neq 0\}}\delta_{(t,\Delta X(t))}$. 
$N$ is a random Poisson measure  on $[0,\infty)\times\{\RL^d\setminus 0\}$ with the mean $ds\times Q(dx)$,
where $Q(dx)$ is the so-called L\'evy measure on $\RL^d\setminus 0$ and $a\in\RL^d$. The abbreviation $N_t(dx)$ means $N([0,t], dx)$.

In this article a stochastic process describing a storage level in storage and inventory systems is investigated. The presented research will distinguish between two cases. The first one will assume that the inflow of goods or products (we do not distinguish products that is we assume that we store one kind of the products) to a storehouse is random due to a L\'evy process and the outflow is deterministic and linear.
Such a situation can happen when the products peel from a conveyor belt, natural resources are mined and then stored or in e.g. fishing industries. A similar situation we can have in terminals, dams or telecommunication systems where the inflow is random and the outflow is deterministic. This case is particularly interesting as we know from the review by Silver \cite{silver08} that it is uncommon in the scientific literature to model supplies (inflow) for storage systems in this manner. 

In the second case we assume conversely that is the inflow is deterministic and the outflow is random.
Such a situation happens when the customers buy a certain product and a storehouse can guarantee deterministic linear supplies and the demands of customers are random. In this case we consider the so-called inventory systems. From Silver \cite{silver08} we know that these assumptions are quite common in research. 

In this paper we want to review theoretical results which determines the probability that the storage level
exceeds a certain value.  Such a characteristic is important in the management of storage and inventory systems. In many cases this value can be determined from computer simulation but to check the accuracy of the simulation we need theoretical formulas. Moreover the theoretical results are interesting from the statistical point of view, that is they are necessary to fit a model to real data. They are also noteworthy from the mathematical point of view as they are based on deep results from the theory of stochastic processes. 

Let us start to describe our models more precisely. In the first case we assume that the supply to a storehouse is random and due to a L\'evy process $X$ that is the whole supply to a storehouse until time $t$ is $X(t)$. From the practical point of view we should assume that the processes $X$ is non-negative with non-decreasing sample paths but we consider several examples where these assumptions will be dropped.
Then we can regard that some uncertainty (a noise) is added to the outflow from the storehouse.
First for the simplicity let us assume that the net supply is the following
$$
Y(t)=X(t)-t\,,
$$
where $X$ is a L\'evy process with values in $[0,\infty)$ and with non-decreasing sample paths. 
Thus the process $Z(t)$ describing the stock level of products in the storehouse at time $t$ satisfies
the following equation
\begin{equation}\label{gr}
Z(t)=Z(0)+Y(t)+\int_0^t \ind\{Z(s)=0\}\,ds\,,
\end{equation}
where $Z(0)$ is the stock level at time $t=0$ (it can be a random or deterministic quantity) and 
$$
 \ind\{Z(s)=0\}=\left\{
\begin{array}{ll}
1 &\mbox{if}\,\, Z(s)=0\,,\\
0 &\mbox{if}\,\, Z(s)>0
\end{array}
\right.
$$
which means that the stock level can not be less than zero and if $Z(s)=0$ the supplies are not greater than demands which are linear. The next proposition gives a solution for eq. (\ref{gr}) (see e.g. Prabhu \cite{pr:98}).
\begin{proposition}\label{ssup}
The integral equation (\ref{gr}) has a unique solution in the following form
\begin{equation}\label{procg}
Z(t)=Z(0)+Y(t)+I(t)\,,
\end{equation}
where
$$
I(t)=[\inf_{s\leq t} Y(s)+Z(0)]^-\,.
$$
Moreover if $Z(0)=0$ then 
\begin{equation}\label{sup}
Z(t)=\sup_{s\leq t}[Y(t)-Y(s)]\stackrel{d}{=}\sup_{s\leq t}Y(s)\,.
\end{equation}
\end{proposition}
\proof
Subtracting by the sides the equation (\ref {gr}) for times $t$ and $\tau-$  (the left limit at $\tau$, $\tau\leq t$) we get
\begin{eqnarray*}
Z(t)&=&Z(\tau-)+Y(t)-Y(\tau-)+\int_\tau^t\ind\{Z(s)=0\}\,ds\\
&\geq& Y(t)-Y(\tau-)\,,
\end{eqnarray*}
where the last inequality follows from the fact that the integral in the equation is non-negative. First let us assume that the storehouse is empty at some moments until time $t$ (precisely on the time interval $(0,t]$). Then we can define the last moment when the storehouse is empty (until the moment $t$) that is $t_0=\max\{\tau: \tau\leq t,\, Z(\tau-)=0\}$. Thus $t_0$ is the last moment when the storehouse was empty on the time interval $(0,t]$ and
$$
Z(t)=Y(t)-Y(t_0-)
$$
using the last equation and the fact that $\int_{t_0}^t\ind\{Z(s)=0\}\,ds=0$. Hence by the last inequality we obtain 
\begin{equation}\label{rg2}
Z(t)=\sup_{\tau\leq t}[Y(t)-Y(\tau-)]
\end{equation}
(we put $Y(0-)=0$).
Now let us assume that the storehouse is never empty on the time interval $(0,t]$. Then by  eq. (\ref{gr}) it follows that $Z(0)+Y(\tau-)>0$ for all $0<\tau\leq t$ hence $Z(t)=Z(0)+Y(t)>Y(t)-Y(\tau-)$.
Thus by (\ref{rg2}) we get
\begin{eqnarray*}
Z(t)&=&\max\{\sup_{\tau\leq t}[Y(t)-Y(\tau-)], Z(0)+Y(t)\}\\
&=&\max\{Y(t)-\inf_{\tau\leq t}Y(\tau), Z(0)+Y(t)\}\,,
\end{eqnarray*}
where in the last equality we use the stochastic continuity of the process $Y$.
Thus using the last equality and eq. (\ref{gr}) we get
\begin{eqnarray*}
\int_0^t \ind\{Z(s)=0\}\,ds&=& Z(t)-Z(0)-Y(t)\\
&=&\max\{-\inf_{\tau\leq t}Y(\tau)-Z(0), 0\}\\
&=&[\inf_{\tau\leq t}Y(\tau)+Z(0)]^-\,.
\end{eqnarray*}
Now let us assume $Z(0)=0$ and notice that $Y_1(\tau)\stackrel{d}{=}Y(\tau)$ for $\tau\leq t$ (in sense of finite dimensional distributions), where $Y_1(\tau)=Y(t)-Y(t-\tau)$. Then for a given $t$ we have
\begin{eqnarray*}
Z(t)&=&Y(t)-\inf_{s\leq t}Y(s)\\
&=&\sup_{s\leq t}[Y(t)-Y(s)]\\
&=&\sup_{s\leq t}Y_1(t-s)\\
&=&\sup_{s\leq t}Y_1(s)\\
&\stackrel{d}{=}&\sup_{s\leq t}Y(s)
\end{eqnarray*}
where the last equality is in the sense of one dimensional distributions. The proof is completed.
\begin{remark}
Except the eq. (\ref{sup}) the statements of the above theorem are valid for any process $X$ which starts from zero and has non-decreasing trajectories and is stochastically continuous that is we can drop the assumption that $X$ is a L\'evy process.
\end{remark}
\begin{remark}
The equality in (\ref{sup}) is not in the sense of finite dimensional distributions but in the sense
of one dimensional distributions.
\end{remark}
\begin{remark}
The process (\ref{procg}) can be used as a storage level process even though the process $X$ has values in
$\RL$ and its sample paths are not non-decreasing. Then we add some uncertainty to the outflow from the storehouse and the process (\ref{procg}) is a solution of a slightly different equation than (\ref{gr}). Thus the process (\ref{procg}) is more general and can serve as a starting point to describe storage levels.
\end{remark}

As we mentioned we will investigate a model where the supplies are described by a linear function and
the customers demands (outflow) are random. More precisely we assume that the net supplies are the following
$$
Y(t)=t-X(t)\,,
$$
where $X$ is a L\'evy process with values in $[0,\infty)$ and then its sample paths are non-decreasing. Hence the process $Z(t)$ describing the stock level at a moment $t$ satisfies the following equation 
\begin{equation}\label{gr1}
Z(t)=Z(0)+Y(t)+\int_0^t \ind\{Z(s)=0\}\,dX(s)\,,
\end{equation}
where $Z(0)$ is the stock level at time $t=0$ (it can be a deterministic or random quantity). 
\begin{proposition}\label{ssup1}
The integral equation (\ref{gr1}) has a unique solution of the form
\begin{equation}\label{procg1}
Z(t)=Z(0)+Y(t)+I(t)\,,
\end{equation}
where
$$
I(t)=[\inf_{s\leq t} Y(s)+Z(0)]^-\,.
$$
Moreover if $Z(0)=0$ then 
\begin{equation}\label{sup1}
Z(t)=\sup_{s\leq t}[Y(t)-Y(s)]\stackrel{d}{=}\sup_{s\leq t}Y(s)\,.
\end{equation}
\end{proposition}
\proof
The proof goes a similar way as the proof of Prop. \ref{ssup}.
\begin{remark}
Except eq. (\ref{sup1}) the statements of the above theorem are valid for any process $X$ which starts from zero and has non-decreasing trajectories and is stochastically continuous, that is we can drop the assumption that $X$ is a L\'evy process.
\end{remark}
\begin{remark}
The equality in (\ref{sup1}) is not in the sense of finite dimensional distributions but in the sense
of one dimensional distributions.
\end{remark}
\begin{remark}
The process (\ref{procg1}) can be used as a storage level process even though the process $X$ has values in
$\RL$ and its sample paths are not non-decreasing. Then we add a some uncertainty to the inflow to the storehouse and the process (\ref{procg1}) is a solution of a slightly different equation than (\ref{gr1}). Thus the process (\ref{procg1}) is more general and can serve as a starting point to describe storage levels.
\end{remark} 

In this paper we are interested in two quantities that is the probability of an overflow at time $t$ which is
$\Prob(Z(t)>u)$ and the probability of an overflow up to time $t$ which is $\Prob(\sup_{s\leq t}Z(s)>u)$
under assumption that the initial storage level is deterministic $Z(0)=z\geq 0$. In several cases we present explicit formulas for those probabilities or we give Laplace transforms of related quantities. 
In the article Dębicki and Mandjes \cite{de:ma:12} a similar treatment is conducted but they consider
the workload process (in our article it is the process $Z$) in the stationarity regime. More precisely
the random variable $\sup_{s<\infty}Y(s)$ is investigated. They review results which present Laplace transforms of the distribution of $\sup_{s<\infty}Y(s)$ for a spectrally positive or negative L\'evy process $Y$. Moreover the results on the distribution of $Z$ stopped at an exponential time
are presented and the correlation of the workload process $Z$ is derived. In the paper of Dębicki and Mandjes \cite{de:ma:12} they also investigate $\inf_{s\leq t}Z(s)$ but as before under assumption that the workload process is stationary. Almost all their results are in the language of Laplace transforms or double Laplace transforms. Thus our work is disjoint with the work of Dębicki and Mandjes \cite{de:ma:12} and we emphasis on explicit formulas.

\section{Probability of an overflow in storage systems}
In this section we will investigate the process $Z$ from Prop. \ref{ssup}. Thus under assumption $Z(0)=0$ (see eq. (\ref{sup})) we are interested in the following probability
$$
\Prob(Z(t)>u)=\Prob(\sup_{s\leq t}(X(s)-s)>u)
$$
which is the probability that the stock level exceeds $u$ at time $t$.
First we assume that the process $X$ is spectrally positive, i.e. its L\'evy measure $Q$ has the support on the positive half line (this means that its jumps are only positive) and has a finite variation which is equivalent for the L\'evy process $X$ that 
$\int_0^1x\,Q(dx)<\infty$. The next theorem is a direct consequence of Michna \cite{mi:11} and \cite{mi:11a} which is originally due to Tak\'acs \cite{ta:65} in a little more general setting  (for the processes with non-negative interchangeable increments that is in our case $X$ has to have the finite variation).
\begin{theorem}\label{finitem} 
Let $u>0$, $Z(0)=0$ and $X$ be a spectrally positive L\'evy process with the finite variation. Then
\begin{eqnarray*}
\Prob(Z(t)>u)&=&\Prob(\sup_{s\leq t}(X(s)-s)>u)\\
&=&\Prob(X(t)-t > u)+\\
&&\int_0^t\frac{f(u+s,s)}{t-s}
\,ds\int_{0}^{t-s}\Prob(X(t-s)\leq x)\,dx\,,
\end{eqnarray*}
where $f(x,s)$ is the density function of $X(s)$ provided that it exists and  
$$
\Prob(Z(t)>0)=1-\frac{1}{t}\int_0^{t}\Prob(X(t)\leq x)\,dx\,.
$$
\end{theorem}

\begin{example}
Let $X$ be gamma process. Then its L\'evy measure is the following
$$
Q(dx)=\frac{a\ind\{x>0\}}{x}\,\exp\left(-\frac{x}{b}\right)\,dx\,,
$$
where $a>0$ and $b>0$. This process is a subordinator e.g. its trajectories are non-decreasing. 
The density function of $X(s)$ is of the form
$$
f(x,s)=
\frac{1}{b^{as}\Gamma(as)}\,x^{as-1}\exp\left(-\frac{x}{b}\right)\ind\{x>0\}\,.
$$
Thus the quantity of the inflow until time $s$ is distributed due to the above density function. 
Hence using Th. \ref{finitem} we get for $u>0$ (see also Dickson and Waters \cite{di:wa:93})
\begin{eqnarray*}
\lefteqn{\Prob(Z(t)> u)}\\
&=&\Prob(X(t)-t>u)\\
&&\,\,\,\,+\int_0^t \Prob(X(t-s)\leq t-s)f(u+s,s)\,ds\\
&&\,\,\,\,-ab\int_0^t\Prob(X(t-s+1/a)\leq t-s)f(u+s,s)\,ds\,,
\end{eqnarray*}
where
$$
\Prob(X(s)\leq x)=\frac{1}{b^{as}\Gamma(as)}\int_0^x y^{as-1}\exp\left(-\frac{y}{b}\right)\,dy\,. 
$$
\end{example}

\begin{example}
Let $X$ be inverse Gaussian L\'evy process. Then its L\'evy measure is the following
$$
Q(dx)=\frac{\delta\ind\{x>0\}}{\sqrt{2\pi x^3}}\exp(-\frac{1}{2}\gamma x)\,dx
$$
where $\delta>0$ and $\gamma>0$. This process is a subordinator.  
The density function of $X(s)$ is the following
$$
f(x,s)=\frac{\delta s e^{\gamma\delta s}}{\sqrt{2\pi}}x^{-3/2}\exp\left\{-\frac{1}{2}(\delta^2 s^2 x^{-1}+\gamma^2 x)\right\}\ind\{x>0\}\,.
$$
As in the previous example for $u>0$ we obtain
\begin{eqnarray*}
\lefteqn{\Prob(Z(t)> u)}\\
&=&\frac{\delta te^{\gamma\delta t}}{\sqrt{2\pi}}\int_{u+t}^{\infty}x^{-3/2}\exp\left\{-\frac{1}{2}(\delta^2 t^2x^{-1}+\gamma^2 x)\right\}\,dx\\
&&+\frac{\delta^2 e^{\gamma\delta t}}{2\pi}\int_0^{t}s(u+s)^{-3/2}\exp\left\{-\frac{1}{2}(\delta^2s^2(u+s)^{-1}+\gamma^2 (u+s))\right\}\,ds\\
&&\,\,\,\,\,\,\cdot\int_0^{t-s}((t-s) x^{-3/2}-x^{-1/2})\exp\left\{-\frac{1}{2}(\delta^2(t-s)^2x^{-1}+\gamma^2 x)\right\}\,dx\,.
\end{eqnarray*}
Similarly one can consider generalized inverse Gaussian L\'evy process.
\end{example}

Now let us investigate the case when $X$ is a spectrally positive L\'evy process with the infinite variation which means that $\int_0^1x\,Q(dx)=\infty$. Let us notice that in this case the trajectories of $X$ are not non-decreasing which means that in our model we add a some uncertainty to the outflow. The next result is due to Michna \cite{mi:11}
and \cite{mi:11a}.
\begin{theorem}\label{infinitem}
Let $u>0$, $Z(0)=0$ and $X$ be a spectrally positive L\'evy process with the infinite variation. Then
\begin{eqnarray*}
\Prob(Z(t)>u)&=&\Prob(\sup_{s\leq t}(X(s)-s)>u)\\
&=&\Prob(X(t)-t > u)+\\
&&\int_0^t\frac{f(u+s,s)}{t-s}
\,ds\int_{-\infty}^{t-s}\Prob(X(t-s)\leq x)\,dx\,,
\end{eqnarray*}
where $f(x,s)$ is the density function of  $X(s)$ provided that it exists and 
$$
\Prob(Z(t)>0)=1\,.
$$
\end{theorem}
For the proof see Michna \cite{mi:11} and \cite{mi:11a}.
\begin{example}
Let $X$ be an $\alpha$-stable spectrally positive L\'evy process with $1\leq \alpha<2$. Its L\'evy measure is of the form
$$
Q(dx)=\frac{\sigma\ind\{x>0\}}{x^{\alpha+1}}\,dx
$$
where $\sigma>0$.
The density function of $X(s)$ exists but it can not be presented as an elementary function (see e.g.
Zolotareav \cite{zo:86} or Nolan \cite{no:97}). There exist integral or series forms of $f(x,s)$ so applying Th. \ref{infinitem} one can compute numerically probability that the stock level will be greater than $u$ at time $t$. 
\end{example}
\begin{example}
Let $X$ be a tempered $\alpha$-stable spectrally positive L\'evy process with $1\leq \alpha<2$ that is its L\'evy measure is the following (see Koponen \cite{ko:95})
$$
Q(dx)=\frac{\sigma e^{-\lambda x}\ind\{x>0\}}{x^{\alpha+1}}\,dx
$$
where $\sigma>0$ and $\lambda>0$.
The density function of $X(s)$ one can determine by the inverse Fourier transform of the characteristic function. Thus by Th. \ref{infinitem} we can compute numerically probability that the stock level will be greater than $u$ at time $t$.
\end{example}

According to Prop. \ref{ssup} the stock level at time $t$ for $Z(0)=z$ ($z\geq 0$ is deterministic) is the following
$$
Z(t)=z+X(t)-t+[\inf_{s\leq t}(X(s)-s)+z]^-\,.
$$
Let us put
\begin{equation}\label{expy}
\Exp e^{-wY(1)}=e^{\psi(w)}\,,
\end{equation}
where $w\geq 0$ and
$$
\tau(u)=\inf\{t\geq 0: Z(t)>u\}\,.
$$
Let us notice that $\tau(u)$ is the first epoch when the stock level exceeds $u$. Moreover probability that the stock level exceeds $u$ until time $t$ is 
$$
\Prob(\sup_{s\leq t}Z(s)>u)=\Prob(\tau(u)< t)\,.
$$
Next theorem is due to Pistorius (see Pistorius \cite{pi:04} see also Avram et al. \cite{av:ky:pi:04})
\begin{theorem}\label{mt}
For $0\leq z\leq u$
$$
\Exp e^{-r\tau(u)}=K^{(r)}(u-z)-rW^{(r)}(u-z)W^{(r)}(u)/{W^{(r)}}'{_{+}}(u)\,,
$$
where $r\geq 0$
and
$$
K^{(r)}(x)=1+r\int_0^x W^{(r)}(y)\,dy\,,
$$
where $W^{(r)}(y)$ is the so-called $r$-scale function on $[0,\infty)$ determined by Laplace transform
$$
\int_0^\infty e^{-wy}W^{(r)}(y)\,dy=\frac{1}{\psi(w)-r}
$$
defined for $w>\Phi(r)$, where $\Phi(\cdot)$ is the right inverse function to $\psi(\cdot)$ (see eq. (\ref{expy})) and ${W^{(r)}_+}'(u)$ is the right derivative of $W^{(r)}(x)$ at $x=u$.
\end{theorem}
The above theorem gives the full solution to the problem of the overflow of the stock level in the model with the deterministic linear outflow. The disadvantage of the above formula is that we can obtain the probability only numerically by inverting Laplace transform of complicated functions.
 
 From the above theorem one can easily find the expected value of the first epoch when the store level exceeds $u$. Let us put $W(x)=W^{(0)}(x)$.
\begin{corollary}
$$
\Exp \tau(u)=W(u-z)\frac{W(u)}{{W'}_{+}(u)}-\bar{W}(u-z)\,.
$$
where ${W'}_{+}(u)$ is the right derivative of $W(x)$ at $x=u$ and $\bar{W}(x)=\int_0^x W(y)\,dy$.
\end{corollary}

\section{Probability of an overflow in inventory systems}
In this section we will assume that the supplies are linear and the outflow is a spectrally positive L\'evy process (its jumps are only positive) that is the net supplies until a moment $t$ are the following
$$
Y(t)=t-X(t)
$$
where $X$ is a spectrally positive L\'evy process. Let us notice that $Y$ is now a spectrally negative 
L\'evy process.
From the mathematical and practical point of view the situation is different than in the model with the linear outflow. Under these assumptions we give a complete solution for any initial stock level. This is highly desirable as it gives us a general model, for any outflow modeled as a L\'evy process. Furthermore this is also interesting as this probability is directly linked to the term of a service level within inventory management see Silver et al. \cite{silveretal98}. The service level corresponds to the percentage of customer demand (in the form of volume) serviced directly from stock. In the simplified case when $Z(0)=0$ the following probability
$$
\Prob(Z(t)>u)=\Prob(\sup_{s\leq t}(s-X(s))>u)
$$
can be derived from Kendall's identity which is valid for spectrally negative L\'evy processes, here $Y(t)=t-X(t)$ is spectrally negative (see Kendall \cite{ke:57} or e.g. Borovkov and Burq
 \cite{bo:bu:01} and the references therein, see also Tak\'acs \cite{ta:65} for the case when $X$ has the finite variation).  
For an elementary proof of an explicit formula for the above probability in the case of the spectrally negative $\alpha$-stable L\'evy process without any linear drift with $1<\alpha\leq 2$
see Michna \cite{mi:13}.
In a more general setting we find the probability that the stock level will exceed $u$ until the moment $t$ that is
$$
\Prob(\sup_{s\leq t}Z(s)>u)
$$
if $Z(0)=z$ for $0\leq z\leq u$ and $z$ is deterministic. However as before the solution will not be in an explicit form but as Laplace transformations of scale functions.

According to Prop. \ref{ssup1} the stock level at time $t$ is the following
$$
Z(t)=z+t-X(t)+[\inf_{s\leq t}(s-X(s))+z]^-\,.
$$
Let us put
\begin{equation}\label{expy1}
\Exp e^{wY(1)}=e^{\psi(w)}\,,
\end{equation}
where $w\geq 0$ and
$$
\tau(u)=\inf\{t\geq 0: Z(t)>u\}\,.
$$
Let us notice that $\tau(u)$ is the first epoch when the stock level exceeds $u$. Moreover probability that the stock level exceeds $u$ until time $t$ is 
$$
\Prob(\sup_{s\leq t}Z(s)>u)=\Prob(\tau(u)< t)\,.
$$
Next theorem is due to Pistorius (see Pistorius \cite{pi:04})
\begin{theorem}\label{mt1} 
For $0\leq z\leq u$
$$
\Exp e^{-r\tau(u)}=\frac{K^{(r)}(z)}{K^{(r)}(u)}\,,
$$
where $r\geq 0$
and
$$
K^{(r)}(x)=1+r\int_0^x W^{(r)}(y)\,dy\,,
$$
where $W^{(r)}(y)$ is the so-called $r$-scale function on $[0,\infty)$ determined by Laplace transform
$$
\int_0^\infty e^{-wy}W^{(r)}(y)\,dy=\frac{1}{\psi(w)-r}
$$
defined for $w>\Phi(r)$, where $\Phi(\cdot)$ is the right inverse function to $\psi(\cdot)$ (see eq. (\ref{expy1})).
\end{theorem}
\begin{remark}
By eq. (\ref{expy}) and (\ref{expy1}) the function $\psi$ and the scale functions in Th. \ref{mt} and
\ref{mt1} are the same.
\end{remark}
Similarly as in Th. \ref{mt} the above formula permits to determine the probability of the overflow of the storehouse but we do not avoid numerical calculations. For some L\'evy processes the scale function is available in an explicit form that is for compound Poisson processes, $\alpha$-stable processes with $\alpha\in (0,2)$, tempered $\alpha$-stable processes and Brownian motion (see e.g. \cite{bi:ky:10} and the references therein). Generally one can not obtain an explicit form of scale functions. However there are numerical methods which permit to determine the scale function of a general spectrally negative L\'evy process (see Surya \cite{su:08} for Laplace inverse algorithm).

From the above theorem one can easily find the expected value of the first epoch when the store level exceeds $u$. 
\begin{corollary}
$$
\Exp \tau(u)=\bar{W}(u)-\bar{W}(z)\,.
$$
where $W(x)=W^{(0)}(x)$ and $\bar{W}(x)=\int_0^x W(y)\,dy$.
\end{corollary}
An interesting aspect of this research is the relative ease with which asymmetrical distributions of the outflow can be incorporated. This is interesting since we know from e.g. Bobko and Whybark \cite{bobkoandwhybark85} and Zotteri \cite{zotteri00} that the asymmetrical distributions are particularly difficult to manage and tend to lead to high costs of the systems.

\end{document}